# On randomized stopping

ISTVÁN GYÖNGY* and DAVID ŠIŠKA**

*School of Mathematics, University of Edinburgh, King's Buildings, EH9 3JZ, U.K.*
*E-mail: *I.Gyongy@ed.ac.uk; **D.Siska@ed.ac.uk*

A general result on the method of randomized stopping is proved. It is applied to optimal stopping of controlled diffusion processes with unbounded coefficients to reduce it to an optimal control problem without stopping. This is motivated by recent results of Krylov on numerical solutions to the Bellman equation.

## 1. Introduction

It is known that optimal stopping problems for controlled diffusion processes can be transformed into optimal control problems by using the method of randomized stopping (see [2] and [8]). Since only a few optimal stopping problems can be solved analytically (see [13]), one has to resort to numerical approximations of the solution. In such cases, one would like to know the rate of convergence of these approximations. Embedding optimal stopping problems into the class of stochastic control problems allows one to apply numerical methods developed for stochastic control [4]. The price one pays for this is the unboundedness of the reward function, as a function of the control parameter.

Recently, a major breakthrough has been made in estimating the rate of convergence of finite difference approximations for the pay-off functions of stochastic control problems (in [9], followed by [10] and [11]). Applying Krylov's methods, new rate of convergence estimates can be found in [1, 5, 6, 7]. New estimates applicable to numerical approximations of normalized Bellman equations appear in [12].

Our main result, Theorem 2.1, formulates the method of randomized stopping in a general setting. Applying it to optimal stopping problems of controlled diffusion processes we easily get (see Theorem 3.2) that under general conditions, the pay-off function of optimal stopping problem of controlled diffusions equals the pay-off function of the control problem obtained by randomized stopping. This result is known from [8] in the case where the coefficients of the controlled diffusions are bounded in the control parameter (see Section 4 of Chapter 3 in [8]). In Theorem 3.2, the coefficients of the diffusions and the functions defining the pay-off may be unbounded functions of the control parameter. Also, they need not satisfy those conditions on continuity which are needed in [8]. Theorem 2.1 can also be applied to optimal stopping of stochastic processes from a larger classes than that of diffusion processes. For the theory of controlled diffusion processes, we refer to [8].

## 2. Main result

Let $(\Omega, \mathcal{F}, \mathbb{P})$ be a probability space equipped with a filtration $(\mathcal{F}_t)_{t\geq 0}$, such that $\mathcal{F}_0$ contains all $\mathbb{P}$-null sets. Let $\mathfrak{T}$ denote the set of finite stopping times. Let $\mathfrak{F}$ denote the set of all processes $F = (F_t)_{t\geq 0}$ which are $\{\mathcal{F}_t\}_{t\geq 0}$ adapted, right-continuous and increasing, such that

$$\forall \omega \in \Omega \qquad F_0(\omega) = 0 \quad \text{and} \quad \lim_{t\to\infty} F_t(\omega) = 1.$$

Let $\bar{\mathfrak{R}}$ be a class of non-negative adapted locally integrable stochastic processes $r = (r_t)_{t\geq 0}$ such that

$$\int_0^\infty r_t \, dt = \infty.$$

Let $\mathfrak{R}_n$ denote those stochastic processes from $\bar{\mathfrak{R}}$ which take values in $[0,n]$. Set $\mathfrak{R} = \bigcup_{n\in\mathbb{N}} \mathfrak{R}_n$.

**Theorem 2.1.** *Let $(h_t)_{t\geq 0}$ be a progressively measurable process with sample paths continuous at 0. Assume that for all $t \geq 0$, $|h_t| \leq \xi$ for some random variable $\xi$ satisfying $\mathbb{E}\xi < \infty$. Then,*

$$\sup_{\tau \in \mathfrak{T}} \mathbb{E} h_\tau = \sup_{F \in \mathfrak{F}} \mathbb{E} \int_0^\infty h_t \, dF_t < \infty. \tag{2.1}$$

Theorem 2.1 is will be proven at the end of this section.

**Theorem 2.2.** *Let $(h_t)_{t\geq 0}$ be an adapted cadlag process such that*

$$\mathbb{E} \sup_{t\geq 0} |h_t| < \infty.$$

*Then,*

$$\sup_{r\in\bar{\mathfrak{R}}} \mathbb{E} \int_0^\infty h_t r_t e^{-\int_0^t r_u du} \, dt = \sup_{r\in\mathfrak{R}} \mathbb{E} \int_0^\infty h_t r_t e^{-\int_0^t r_u du} \, dt = \sup_{\tau\in\mathfrak{T}} \mathbb{E} h_\tau < \infty. \tag{2.2}$$

**Proof.** Let $r \in \bar{\mathfrak{R}}$. If $\varphi_t = e^{-\int_0^t r_u du}$, then $1 - \varphi \in \mathfrak{F}$. Hence

$$\mathbb{E} \int_0^\infty h_t r_t e^{-\int_0^t r_u du} \, dt \leq \sup_{F\in\mathfrak{F}} \mathbb{E} \int_0^\infty h_t \, dF_t = \sup_{\tau\in\mathfrak{T}} \mathbb{E} h_\tau < \infty,$$

by Theorem 2.1. On the other hand, for $\tau \in \mathfrak{T}$, let $r^n_t = 0$ for $t < \tau$ and $r^n_t = n$ for $t \geq \tau$. Set

$$F^n_t := \begin{cases} 0, & \text{for } t < \tau, \\ 1 - e^{-n(t-\tau)}, & \text{for } t \geq \tau. \end{cases}$$

Then, for any $\omega \in \Omega$ and any $\delta > 0$,

$$\int_0^\infty h_t r^n_t e^{-\int_0^t r^n_u du} \, dt = \int_0^\infty h_t \, dF^n_t = \int_\tau^\infty h_t \, dF^n_t = I_n + J_n + K_n,$$

where

$$I_n := \int_\tau^{\tau+\delta} h_\tau \, dF^n_t, \qquad J_n := \int_\tau^{\tau+\delta} (h_t - h_\tau) \, dF^n_t, \qquad K_n := \int_{\tau+\delta}^\infty h_t \, dF^n_t.$$

Notice that as $n \to \infty$,

$$I_n = h_\tau (1 - e^{-n\delta}) \to h_\tau \quad \text{and} \quad |K_n| \leq \sup_{t \geq 0} |h_t| e^{-(\tau+\delta)n} \to 0.$$

Furthermore,

$$|J_n| \leq \sup_{t \in [\tau, \tau+\delta]} |h_t - h_\tau| \to 0 \quad \text{as } \delta \to 0.$$

Hence, for any $\omega \in \Omega$, first taking the limit as $n \to \infty$ and then the limit as $\delta \to 0$,

$$\int_0^\infty h_t r^n_t e^{-\int_0^t r^n_u du} \, dt \to h_\tau.$$

Clearly,

$$\left| \int_0^\infty h_t r^n_t e^{-\int_0^t r^n_u du} \, dt \right| \leq \sup_{t \geq 0} |h_t|.$$

Hence, by Lebesgue's theorem,

$$\sup_{r \in \mathfrak{R}} \mathbb{E} \int_0^\infty h_t r_t e^{-\int_0^t r_u du} \, dt \geq \lim_{n \to \infty} \mathbb{E} \int_0^\infty h_t r^n_t e^{-\int_0^t r^n_u du} \, dt$$

$$= \mathbb{E} \lim_{n \to \infty} \int_0^\infty h_t r^n_t e^{-\int_0^t r^n_u du} \, dt = \mathbb{E} h_\tau. \qquad \square$$

**Proof of Theorem 2.1.** Due to the assumptions on $h_t$ in Theorem 2.1,

$$\sup_{F \in \mathfrak{F}} \mathbb{E} \int_0^\infty h_t \, dF_t \leq \sup_{F \in \mathfrak{F}} \mathbb{E} \int_0^\infty \xi \, dF_t = \mathbb{E} \xi < \infty.$$

Let $\tau \in \mathfrak{T}$. Define

$$G^n_t := \mathbf{1}_{\{\tau > 0\}} \mathbf{1}_{\{t \geq \tau\}} + \mathbf{1}_{\{\tau = 0\}} (1 - e^{-nt}).$$

Clearly, $G_t^n \in \mathfrak{F}$. Using the argument of the proof of Theorem 2.2, it is easy to see that

$$\mathbb{E}h_\tau = \lim_{n\to\infty} \mathbb{E}\int_0^\infty h_t \, dG_t^n \leq \sup_{F\in\mathfrak{F}} \mathbb{E}\int_0^\infty h_t \, dF_t.$$

To complete the proof, one needs to show that

$$\sup_{\tau\in\mathfrak{T}} \mathbb{E}h_\tau \geq \mathbb{E}\int_0^\infty h_t \, dF_t$$

holds for any $F \in \mathfrak{F}$. To this end, let us define the time-change $\beta$ by

$$\beta(r) := \inf\{t \geq 0 : F_t \geq r\}.$$

It is then easy to see that $\beta(r)$ is a stopping time for every $r \geq 0$ and that

$$M := \mathbb{E}\int_0^\infty h_t \, dF_t = \mathbb{E}\int_0^1 h_{\beta(r)} \, dr = \int_0^1 \mathbb{E}h_{\beta(r)} \, dr$$
$$\leq \mathbb{E}\int_0^1 \sup_{\tau\in\mathfrak{T}} \mathbb{E}h_\tau \, dr = \sup_{\tau\in\mathfrak{T}} \mathbb{E}h_\tau. \qquad \Box$$

(2.3)

**Remark 2.3.** Notice that from (2.3), we have the existence of $\bar{r} \in [0,1)$ such that

$$M = \int_0^1 \mathbb{E}h_{\beta(r)} \, dr \leq \mathbb{E}h_{\beta(\bar{r})} = \mathbb{E}h_{\bar{\tau}}$$

for the stopping time $\bar{\tau} = \beta(\bar{r})$. If $h$ is a piecewise constant process with finitely many jumps, then the following lemma provides a way of constructing a stopping time $\tau$ such that $M \leq \mathbb{E}h_\tau$ holds.

**Lemma 2.4.** *Let $(h_i)_{i=1}^n$, $(p_i)_{i=1}^n$ be sequences of random variables adapted to a filtration $(\mathcal{G}_i)_{i=1}^n$, where $\mathcal{G}_1$ contains all $\mathbb{P}$-null sets, such that for all $i$, $p_i \geq 0$, $\mathbb{E}|h_i| < \infty$ and*

$$\sum_{i=1}^n p_i = 1 \qquad a.s.$$

*There then exist disjoint sets $(A_i)_{i=1}^n$, $A_i \in \mathcal{G}_i$, $A_1 \cup A_2 \cup \cdots \cup A_n = \Omega$ such that almost surely*

$$\mathbb{E}(p_1 h_1 + p_2 h_2 + \cdots + p_n h_n | \mathcal{G}_1) \leq \mathbb{E}(h_1 \mathbf{1}_{A_1} + h_2 \mathbf{1}_{A_2} + \cdots + h_n \mathbf{1}_{A_n} | \mathcal{G}_1).$$

**Proof.** For $n = 1$, the statement of the lemma is obvious. Assume $n = 2$ (this will illustrate the general case better). Then, since $p_1 h_1$ is $\mathcal{G}_1$ measurable, $p_2 = 1 - p_1$ is also $\mathcal{G}_1$ measurable and

$$I := p_1 h_1 + \mathbb{E}(p_2 h_2 | \mathcal{G}_1) = p_1 h_1 + p_2 \mathbb{E}(h_2 | \mathcal{G}_1).$$

Let $A_1 = \{h_1 \geq \mathbb{E}(h_2|\mathcal{G}_1)\}$. This is a $\mathcal{G}_1$ set and

$$I \leq h_1 \mathbf{1}_{A_1} + \mathbf{1}_{\Omega \setminus A_1} \mathbb{E}(h_2|\mathcal{G}_1) = \mathbb{E}(h_1 \mathbf{1}_{A_1} + h_2 \mathbf{1}_{\Omega \setminus A_1}|\mathcal{G}_1).$$

Assume that the lemma holds for $n-1 \geq 1$. Let us prove that it remains true for $n$. Let $B = \{p_1 < 1\}$.

$$I := \mathbb{E}(h_1 p_1 + h_2 p_2 + \cdots + h_n p_n|\mathcal{G}_1)$$
$$= \mathbb{E}(\mathbf{1}_B(h_1 p_1 + h_2 p_2 + \cdots + h_n p_n)|\mathcal{G}_1) + \mathbb{E}(\mathbf{1}_{B^c} h_1|\mathcal{G}_1) =: I_1 + I_2.$$

Then,

$$I_1 = \mathbf{1}_B(h_1 p_1 + \mathbb{E}(h_2 p_s + \cdots + h_n p_n|\mathcal{G}_1))$$
$$= \mathbf{1}_B\left(h_1 p_1 + (p_2 + \cdots + p_n)\mathbb{E}\left(\frac{h_2 p_2 + \cdots + h_n p_n}{p_2 + \cdots + p_n} \,\bigg|\, \mathcal{G}_1\right)\right),$$

since $p_2 + \cdots + p_n > 0$ is $\mathcal{G}_1$ measurable. Let

$$A_1 = \left\{h_1 \geq \mathbb{E}\left(\frac{h_2 p_2 + \cdots + h_n p_n}{p_2 + \cdots + p_n} \,\bigg|\, \mathcal{G}_1\right)\right\} \cap B.$$

Since $A_1$ is a $\mathcal{G}_1$ set,

$$I_1 \leq h_1 \mathbf{1}_{A_1} + \mathbf{1}_B \mathbb{E}\left(\mathbf{1}_{\Omega \setminus A_1} \frac{h_2 p_2 + \cdots + h_n p_n}{p_2 + \cdots + p_n} \,\bigg|\, \mathcal{G}_1\right)$$
$$= h_1 \mathbf{1}_{A_1} + \mathbb{E}(h'_2 p'_2 + \cdots + h'_n p'_n|\mathcal{G}_1),$$

where $h'_i = h_i \mathbf{1}_B \mathbf{1}_{\Omega \setminus A_1}$, $p'_i = \frac{p_i}{p_2 + \cdots + p_n}$ on $B$ and $p'_i = (n-1)^{-1}$ on $B^c$, for $2 \leq i \leq n$. Then, $p'_2 + \cdots + p'_n = 1$. Apply the inductive hypothesis to $(h'_i)_{i=2}^n$, $(p'_i)_{i=2}^n$, $(\mathcal{G}_i)_{i=2}^n$. There are then disjoint $A'_2 \cup A'_3 \cup \cdots \cup A'_n = \Omega$ such that $A'_i \in \mathcal{G}_i$ for $i = 2, \ldots, n$ and

$$\mathbb{E}(h'_2 p'_2 + \cdots + h'_n p'_n|\mathcal{G}_1) \leq \mathbb{E}(h'_2 \mathbf{1}_{A'_2} + \cdots + h'_n \mathbf{1}_{A'_n}|\mathcal{G}_1).$$

Hence,

$$I_1 \leq h_1 \mathbf{1}_{A_1} + \mathbf{1}_B \mathbb{E}(\mathbf{1}_{\Omega \setminus A_1}(h_2 \mathbf{1}_{A'_2} + \cdots + h_n \mathbf{1}_{A'_n})|\mathcal{G}_1).$$

We see that $(\Omega \setminus A_1) \cap (A'_2 \cup \cdots \cup A'_n) = \Omega \setminus A_1$. For $1 < i \leq n$, define $A_i = B \cap A'_i \cap (\Omega \setminus A_1)$. Such $A_i$ are disjoint, $\mathcal{G}_i$ measurable and

$$B^c \cup A_1 \cup A_2 \cup \cdots \cup A_n = \Omega.$$

Thus,

$$I_1 \leq \mathbb{E}(h_1 \mathbf{1}_{A_1} + h_2 \mathbf{1}_{A_2} + \cdots + h_n \mathbf{1}_{A_n}|\mathcal{G}_1).$$

Finally,
$$I = I_1 + I_2 \leq \mathbb{E}(h_1(\mathbf{1}_{A_1} + \mathbf{1}_{B^c}) + h_2\mathbf{1}_{A_2} + \cdots + h_n\mathbf{1}_{A_n}|\mathcal{G}_1). \qquad \square$$

Let $h(t) = \sum_{i=1}^{N} \mathbf{1}_{(\tau_{i-1},\tau_i]}(t)\xi_i$ for an increasing sequence of finite random variables $0 = \tau_0 \leq \tau_1 \leq \cdots \leq \tau_N$ and random variables $\xi_i$ such that $\mathbb{E}|\xi_i| < \infty$ for all $i = 1, 2, \ldots, N$. Let $(\mathcal{F}_t)_{t\geq 0}$ be the filtration generated by the process $h$ and the processes $\mathbf{1}_{(0,\tau_i]}$ for $i = 1, 2, \ldots, N$. Then, by the above lemma applied to $h_i := \xi_i$, $p_i := F_{\tau_i} - F_{\tau_{i-1}}$ and $\mathcal{G}_i := \mathcal{F}_{\tau_i}$,

$$\mathbb{E}\int_0^\infty h(t)\,\mathrm{d}F_t = \mathbb{E}\sum_{i=1}^{N} h_i p_i \leq \mathbb{E}\sum_{i=1}^{N} h_i \mathbf{1}_{A_i} = \mathbb{E}h(\tau),$$

for any $F \in \mathfrak{F}$, where $\tau = \sum_{i=1}^{N} \tau_i \mathbf{1}_{A_i}$, with $A_i$ given by the lemma, is a finite stopping time. Hence, we note that the inequality

$$\mathbb{E}\int_0^\infty h_t\,\mathrm{d}F_t \leq \sup_{\tau\in\mathfrak{T}} \mathbb{E}h_\tau$$

follows for the general case of progressively measurable processes $h$ satisfying the condition of Theorem 2.1, by a simple application of the following lemma from [3].

**Lemma 2.5.** *Let $(h_t)_{t\geq 0}$ be a $\mathcal{B}([0,\infty)) \times \mathcal{F}$-measurable process and let $F \in \mathfrak{F}$ such that*

$$\mathbb{E}\int_0^\infty |h_t|\,\mathrm{d}F_t < \infty.$$

*Then, for each integer $n \geq 1$, there exists a finite sequence of stopping times $\tau_i^{(n)}$ such that*

$$0 = \tau_0^{(n)} < \tau_1^{(n)} \leq \cdots \leq \tau_{N(n)}^{(n)} < \infty$$

*and for*

$$h_t^{(n)} = \sum_{i=1}^{N(n)} h_{\tau_i^{(n)}} \mathbf{1}_{(\tau_{i-1}^{(n)},\tau_i^{(n)}]}(t), \qquad (2.4)$$

*one has*

$$\mathbb{E}\int_0^\infty |h_t - h_t^{(n)}|\,\mathrm{d}F_t \to 0 \qquad as\ n \to \infty.$$

## 3. Application to controlled diffusion processes

We first introduce some notions and notations from the theory of controlled diffusion processes from [8]. Fix $T \in [0,\infty)$. Let $(w_t, \mathcal{F}_t)$ be a $d'$-dimensional Wiener martingale, that is, $w = (w_t)_{t\geq 0}$ is a Wiener process and $w_t - w_s$ is independent of $\mathcal{F}_s$ for all $0 \leq s < t$.

Let $A$ be a separable metric space. For every $t \in [0,T]$, $x \in \mathbb{R}^d$ and $\alpha \in A$, we are given a $d \times d'$-dimensional matrix $\sigma^\alpha(t,x)$, a $d$-dimensional vector $\beta^\alpha(t,x)$ and real numbers $c^\alpha(t,x)$, $f^\alpha(t,x)$ and $g(t,x)$.

**Assumption 3.1.** *$\sigma, \beta, c, f$ are Borel functions of $(\alpha, t, x)$. The function $g$ is continuous in $(t,x)$. There exist an increasing sequence of subsets $A_n$ of $A$ and positive real constants $K$, $K_n$ and $m$, $m_n$ such that $\bigcup_{n \in \mathbb{N}} A_n = A$ and for each $n \in \mathbb{N}$, $\alpha \in A_n$,*

$$|\sigma^\alpha(t,x) - \sigma^\alpha(t,y)| + |\beta^\alpha(t,x) - \beta^\alpha(t,y)| \leq K_n|x-y|,$$
$$|\sigma^\alpha(t,x)| + |\beta^\alpha(t,x)| \leq K_n(1+|x|), \quad (3.1)$$
$$|c^\alpha(t,x)| + |f^\alpha(t,x)| \leq K_n(1+|x|)^{m_n}, \qquad |g(t,x)| \leq K(1+|x|)^m$$

*for all $x \in \mathbb{R}^d$ and $t \in [0,T]$.*

We say that $\alpha \in \mathfrak{A}_n$ if $\alpha = (\alpha_t)_{t \geq 0}$ is a progressively measurable process with values in $A_n$. Let $\mathfrak{A} = \bigcup_{n \in \mathbb{N}} \mathfrak{A}_n$. Then, under Assumption 3.1, it is well known that for each $s \in [0,T]$, $x \in \mathbb{R}^d$ and $\alpha \in \mathfrak{A}$, there is a unique solution $\{x_t : t \in [0, T-s]\}$ of

$$x_t = x + \int_0^t \sigma^{\alpha_u}(s+u, x_u) \, dw_u + \int_0^t \beta^{\alpha_u}(s+u, x_u) \, du, \quad (3.2)$$

denoted by $x_t^{\alpha,s,x}$. For $s \in [0,T]$, we use the notation $\mathfrak{T}(T-s)$ for the set of stopping times $\tau \leq T - s$. Define

$$w(s,x) = \sup_{\alpha \in \mathfrak{A}} \sup_{\tau \in \mathfrak{T}(T-s)} v^{\alpha,\tau}(s,x),$$

where

$$v^{\alpha,\tau}(s,x) = \mathbb{E}_{s,x}^\alpha \left[ \int_0^\tau f^{\alpha_t}(s+t, x_t) e^{-\varphi_t} \, dt + g(s+\tau, x_\tau) e^{-\varphi_\tau} \right],$$

$$\varphi_t = \int_0^t c^{\alpha_r}(s+r, x_r) \, dr$$

and $\mathbb{E}_{s,x}^\alpha$ means expectation of the expression following it, with $x_t^{\alpha,s,x}$ in place of $x_t$ everywhere. It is worth noticing that for

$$w_n(s,x) := \sup_{\alpha \in \mathfrak{A}_n} \sup_{\tau \in \mathfrak{T}(T-s)} v^{\alpha,\tau}(s,x),$$

we have

$$w_n(s,x) \uparrow w(s,x).$$

By Theorem 3.1.8 in [8], $w_n(s,x)$ is bounded from above and below. Hence, $w(s,x)$ is bounded from below. However, it can be equal to $+\infty$.

Let $\mathfrak{R}_n$ contain all progressively measurable, locally integrable processes $r = (r_t)_{(t \geq 0)}$ taking values in $[0, n]$ such that $\int_0^\infty r_t \, dt = \infty$. Let $\mathfrak{R} = \bigcup_{n \in \mathbb{N}} \mathfrak{R}_n$.

Next, we prove a theorem which is known from [8] in the special case when $A = A_n$, $K = K_n$, $m = m_n$ for $n \geq 1$ (see Exercise 3.4.12, via Lemma 3.4.3(b) and Lemma 3.4.5(c)). Our proof is a straightforward application of Theorem 2.2. Since Theorem 2.2 is removed from the theory of controlled diffusion processes developed in [8], we do not require that $\sigma$, $b$, $c$, $f$ and $g$ be continuous in $(\alpha, x)$ and continuous in $x$, uniformly in $\alpha$ for each $t$ (conditions which are needed in [8]).

**Theorem 3.2.** *Let Assumption 3.1 hold. Then, for all $(s, x) \in [0, T] \times \mathbb{R}^d$, either both $w(s, x)$ and*

$$\sup_{\alpha \in \mathfrak{A}} \sup_{r \in \mathfrak{R}} \mathbb{E}^\alpha_{s,x} \left\{ \int_0^{T-s} [f^{\alpha_t}(s+t, x_t) e^{-\varphi_t} + r_t g(s+t, x_t) e^{-\varphi_t}] e^{-\int_0^t r_u \, du} \, dt \right.$$
$$\left. + g(T, x_{T-s}) e^{-\varphi_{T-s} - \int_0^{T-s} r_u \, du} \right\}$$

*are finite and equal, or they are both infinite.*

**Proof.** Without loss of generality, we may assume that $s = 0$. Let $r \in \mathfrak{R}$. For $t > T$, let $f^\alpha(t, x) = 0$, $c^\alpha(t, x) = 0$. For fixed $(\alpha_t) \in \mathfrak{A}$, set

$$f_t = f^{\alpha_t}(t, x_t) e^{-\varphi_t}, \qquad \text{for } t \geq 0,$$

$$g_t = \begin{cases} g(t, x_t) e^{-\varphi_t}, & \text{for } t \leq T, \\ g(T, x_T) e^{-\varphi_T}, & \text{for } t > T. \end{cases}$$

Recall that if $r \in \mathfrak{R}$, then $\int_0^\infty r_t \, dt = \infty$. Clearly, $f_t = 0$ for $t > T$, and

$$\int_T^\infty g_t r_t e^{-\int_0^t r_u \, du} \, dt = g_T e^{-\int_0^T r_u \, du}.$$

Thus,

$$\int_0^T (f_t + g_t r_t) e^{-\int_0^t r_u \, du} \, dt + g_T e^{-\int_0^T r_u \, du}$$
$$= \int_0^\infty (f_t + g_t r_t) e^{-\int_0^t r_u \, du} \, dt$$
$$= \int_0^\infty \left( \int_0^t f_s \, ds + g_t \right) r_t e^{-\int_0^t r_u \, du} \, dt,$$

where the last equality comes from integrating by parts. We check that for

$$h_t := \int_0^t f_s \, ds + g_t,$$

for each $\alpha \in \mathfrak{A}$, $\mathbb{E}\sup_{t\geq 0}|h_t| < \infty$ holds. Indeed, if $\alpha \in \mathfrak{A}_n$, then

$$\mathbb{E}\sup_{t\geq 0}|h_t| \leq \mathbb{E}\sup_{t\in[0,T]}\left(\int_0^t |f^{\alpha_s}(s,x_s)|\mathrm{e}^{-\varphi_s}\,\mathrm{d}s + |g(t,x_t)|\mathrm{e}^{-\varphi_t}\right)$$
$$\leq TK_n\mathbb{E}\sup_{t\in[0,T]}(1+|x_t|)^{m_n} + TK\mathbb{E}\sup_{t\in[0,T]}(1+|x_t|)^m < \infty,$$

due to estimates of moments of solutions to SDEs (Theorem 2.5.9 in [8]). Since $g = g(t,x)$ is continuous, $h_t$ is continuous. Hence, by Theorem 2.2,

$$\sup_{r\in\mathfrak{R}}\mathbb{E}\int_0^\infty \left(\int_0^t f_s\,\mathrm{d}s + g_t\right)r_t\mathrm{e}^{-\int_0^t r_u\mathrm{d}u}\,\mathrm{d}t = \sup_{\tau\in\mathfrak{T}}\mathbb{E}\int_0^\tau f_s\,\mathrm{d}s + g_\tau$$
$$= \sup_{\tau\in\mathfrak{T}(T)}\mathbb{E}\int_0^\tau f_s\,\mathrm{d}s + g_\tau$$

because $f_t = 0$ and $g_t = g_T$ for $t > T$, so nothing can be gained or lost by stopping later. Hence, for any $\alpha \in \mathfrak{A}$,

$$\sup_{\tau\in\mathfrak{T}(T)}\mathbb{E}\left\{\int_0^\tau f^{\alpha_t}(t,x_t)\mathrm{e}^{-\varphi_t}\,\mathrm{d}t + g(\tau,x_\tau)\mathrm{e}^{-\varphi_\tau}\right\}$$
$$= \sup_{r\in\mathfrak{R}}\mathbb{E}\left\{\int_0^T [f^{\alpha_t}(t,x_t)\mathrm{e}^{-\varphi_t} + g(t,x_t)r_t]\mathrm{e}^{-\int_0^t r_u\mathrm{d}u}\,\mathrm{d}t\right.$$
$$\left. + g(T,x_T)\mathrm{e}^{-\varphi_T - \int_0^T r_u\mathrm{d}u}\right\},$$

which proves the theorem. $\square$

## Acknowledgements

The authors are grateful to Nicolai Krylov in Minnesota and Peter Bank in Berlin for useful information on the subject of this paper and to the referee for his insightful suggestions.